\documentclass{article}
\usepackage{fullpage}
\usepackage{graphics,psfrag,epsfig}
\usepackage{amsfonts,latexsym,eucal,amsmath,amsthm,amssymb}

\newcommand{\ZZ}{\mathbb{Z}}

\newcommand{\kb}{\bar{k}}

\title{Trapping Regions for the Navier-Stokes Equations}
\author{Craig Alan Feinstein}
\newcommand{\pd}{\partial}
\date{October 25, 2014}
\begin{document}

\maketitle
\begin{scriptsize}
\begin{center} \bigskip \noindent \textbf{Baltimore, Maryland
U.S.A. email: cafeinst@msn.com, BS"D}
\end{center}

\bigskip \noindent \textbf{Abstract:} In 1999, J.C. Mattingly and Ya. G. 
Sinai used elementary methods to prove the existence and uniqueness of 
smooth solutions to the 2D Navier-Stokes equations with periodic boundary 
conditions. And they were almost successful in proving the existence and 
uniqueness of smooth solutions to the 3D Navier-Stokes equations using the 
same strategy. In this paper, we modify their technique to obtain a simpler proof of 
one of their results. We also argue that there is no logical reason why 
the 3D Navier-Stokes equations must always have solutions, 
even when the initial velocity vector field is smooth; if they do always have 
solutions, it is due to probability and not logic.

\bigskip \noindent \textbf{Disclaimer:} This article was authored by Craig Alan Feinstein in his private
capacity. No official support or endorsement by the U.S. Government is intended or should be inferred.

\end{scriptsize}
\clearpage

In this paper, we examine the three-dimensional Navier-Stokes equations, which model the flow of
incompressible fluids:

\begin{align}
  \label{eq:3DNS}
  \frac{\pd u_i}{\pd t} + \sum_{j=1,2,3} u_j
  \frac{\pd u_i}{\pd x_j} &= \nu \Delta u_i
  -\frac{\pd p}{\pd x_i} \qquad i=1,2,3
  \notag \\
  \sum_{i=1,2,3}  \frac{\pd u_i}{\pd x_i} &= 0 \ ,
\end{align}
where $\nu > 0$ is viscosity, $p$ is pressure, $u$ is velocity, and $t>0$ is time. We shall assume that
both $u$ and $p$ are periodic in $x$. For simplicity, we take the period to be one. The first equation is
Newton's Second Law, force equals mass times acceleration, and the second equation is the assumption that
the fluid is incompressible.

Mattingly and Sinai \cite{b:MatSi99} attempted to show that smooth solutions to 3D Navier Stokes
equations exist for all initial conditions $u(x,0)=u^0(x) \in C^{\infty}$ by dealing with an equivalent
form of the Navier-Stokes equations for periodic boundary conditions:
\begin{align}
  \label{eq:3DvortEqn}
  \frac{\partial \omega_i}{\partial t} + \sum_{j=1,2,3} u_j \frac{\partial
  \omega_i}{\partial x_j} =  \sum_{j=1,2,3} \omega_j \frac{\partial
  u_i}{\partial x_j}+ \nu \Delta \omega_i  \qquad i=1,2,3,
\end{align}
where the vorticity $\omega(x,t)=(\frac{\partial u_2}{\partial x_3} -\frac{\partial u_3}{\partial x_2},
\frac{\partial u_3}{\partial x_1} -\frac{\partial u_1}{\partial
  x_3},\frac{\partial u_1}{\partial x_2} -\frac{\partial u_2}{\partial
  x_1})$.

Their strategy was as follows: Represent the equations \eqref{eq:3DvortEqn} as a Galerkin system in
Fourier space with a basis $\{e^{2 \pi i kx}\}_{k \in \ZZ^3}$. A finite dimensional approximation of this
Galerkin system can be associated to any finite subset $\mathcal{Z}$ of $\ZZ^3$ by setting
$u^{(k)}(t)=\omega^{(k)}(t)=0$ for all $k$ outside of $\mathcal{Z}$. For each finite dimensional
approximation of this Galerkin system, consider the system of coupled ODEs for the Fourier coefficients.
Then construct a subset $\Omega(K)$ of the phase space (the set of possible configurations of the Fourier
modes) so that all points in $\Omega(K)$ possess the desired decay properties. In addition, construct
$\Omega(K)$ so that it contains the initial data. Then show that the dynamics never cause the sequence of
Fourier modes to leave the subset $\Omega(K)$ by showing that the vector field on the boundary of
$\Omega(K)$ points into the interior of $\Omega(K)$.

Unfortunately, their strategy only worked for the 3D Navier-Stokes equations when the Laplacian operator
$\Delta$ in \eqref{eq:3DvortEqn} was replaced by another similar linear operator. (Their strategy was in
fact successful for the 2D Navier-Stokes equations.) In this paper, we attempt to apply their strategy to the
original equations \eqref{eq:3DNS}:

Moving to Fourier space where
\begin{align}
  u_i(x,t)=\sum_{k \in \mathcal Z} u_i^{(k)}(t) e^{2 \pi i kx}
 \ \mbox{, } \
  p(x,t)=\sum_{k \in \mathcal Z} p^{(k)}(t) e^{2 \pi i kx}
 \ \mbox{, and } \
  |k|=\sqrt{\sum_{j=1,2,3} k_j^2} \ ,
\end{align}
let us consider the system of coupled ODEs for a finite-dimensional approximation to the Galerkin-system
corresponding to \eqref{eq:3DNS},
\begin{align}
  \label{eq:first}
  \frac{du_i^{(k)}}{dt}=  \Bigl(\sum_{\substack{q+r=k\\
        q, r \in \mathcal{Z}}}
      \sum_{j=1,2,3} - 2 \pi i q_j u_i^{(q)}u_j^{(r)}\Bigr)
  - 4 \pi^2 \nu |k|^2 u_i^{(k)} - 2 \pi i k_i p^{(k)} \qquad
  i=1,2,3
\end{align}
\begin{align}
  \label{eq:fdivfirst}
  \sum_{i=1,2,3} k_i u_i^{(k)}=0,
\end{align}
where $\mathcal Z$ is a finite subset of $\ZZ^3$ in which $u^{(k)}(t)=p^{(k)}(t)=0$ for each $k \in
\ZZ^3$ outside of $\mathcal Z$. Like the Mattingly and Sinai paper, in this paper, we consider a
generalization of this Galerkin-system:
\begin{align}
  \label{eq:finite}
  \frac{du_i^{(k)}}{dt}=  \Bigl(\sum_{\substack{q+r=k\\
        q, r \in \mathcal{Z}}}
      \sum_{j=1,2,3} - 2 \pi i q_j u_i^{(q)}u_j^{(r)}\Bigr)
  - 4 \pi^2 \nu |k|^\alpha u_i^{(k)} - 2 \pi i k_i p^{(k)} \qquad
  i=1,2,3
\end{align}
\begin{align}
  \label{eq:fdiv}
  \sum_{i=1,2,3} k_i u_i^{(k)}=0,
\end{align}
where $\alpha \geq 2$. Multiplying each of the first three equations by $k_i$ for $i=1,2,3$ and adding
the resulting equations together, we obtain
\begin{align}
   \sum_{\substack{q+r=k\\q, r \in \mathcal{Z}}}
   \sum_{\substack{j=1,2,3\\l=1,2,3}}- 2 \pi i k_l q_j
   u_l^{(q)}u_j^{(r)} = 2 \pi i |k|^2 p^{(k)},
\end{align}
since $\sum_{i=1,2,3} k_i \frac{du_i^{(k)}}{dt}=0$ (by equation \eqref{eq:fdiv}). Then substituting the
above calculated expression for $p^{(k)}$ in terms of $u$ into \eqref{eq:finite} we obtain
\begin{align}
\frac{du_i^{(k)}}{dt}=\Bigl[\sum_{\substack{q+r=k\\
        q, r \in \mathcal{Z}}}\sum_{\substack{j=1,2,3\\l=1,2,3}}2 \pi i
 (\frac{k_i k_l}{|k|^2}-\delta_{il})
   q_j u_l^{(q)}u_j^{(r)}\Bigr]  - 4 \pi^2 \nu |k|^\alpha
   u_i^{(k)}\qquad i=1,2,3.
\end{align}
And since $\sum_{j=1,2,3} r_j u_j^{(r)}=0$ and $q_j+r_j=k_j$, we can substitute $k_j$ for $q_j$:
\begin{align}
\label{eq:final} \frac{du_i^{(k)}}{dt}=\Bigl[\sum_{\substack{q+r=k\\
        q, r \in \mathcal{Z}}}\sum_{\substack{j=1,2,3\\l=1,2,3}}2 \pi i
 (\frac{k_i k_l}{|k|^2}-\delta_{il})
   k_j u_l^{(q)}u_j^{(r)}\Bigr]  - 4 \pi^2 \nu |k|^\alpha
   u_i^{(k)}\qquad i=1,2,3.
\end{align}
Now, we state and prove the following theorem:

\bigskip \noindent \textbf{Theorem:}
  \textit{ Let $\{u^{(k)}(t)\}$ satisfy (\ref{eq:final}), where $\alpha>2.5$. And let $1.5<s<\alpha-1$.
Suppose there exists a constant $C_0>0$ such that $|u^{(k)}(0)| \leq C_0 |k|^{-s}$,
for all $k \in \mathbb{Z}^3$. Then there exists a constant $C>C_0$ such that
$|u^{(k)}(t)| \leq C |k|^{-s}$, for all $k \in \mathbb{Z}^3$ and all $t>0$. (The constants, 
$C_0$ and $C$, are independent of the set $\mathcal{Z}$ defining the Galerkin approximation.)}

\bigskip \noindent \textit{Proof:}
By the basic energy estimate (see \cite{b:CoFo88,b:DoGi95,b:Temam79}), there exists a constant $E \geq 0$
such that for each $t \geq 0$ and for any finite-dimensional Galerkin approximation defined by
$\mathcal{Z} \subset \ZZ^3$, we have $\sum_{k \in \mathcal{Z}}\sum_{i=1,2,3} |u_i^{(k)}(t)|^2 \leq E$.
Hence, for any $K>0$, we can find a $C>C_0$ such that $|\Re(u^{(k)})| \leq C|k|^{-s}$ and
$|\Im(u^{(k)})| \leq C|k|^{-s}$, for all $t \geq 0$ and $k \in \ZZ^3$ with $|k| \leq K$. Now let us consider the set,
\begin{align}
  \label{eq:set}
  \Omega(K)=\Bigl\{\Bigl(\Re(u^{(k)}),\Im(u^{(k)})\Bigr)_{k \in \ZZ^3}
  \mbox{ : }|k| > K \mbox{, }
  |\Re(u^{(k)})| \leq C|k|^{-s} \mbox{, and } |\Im(u^{(k)})| \leq
  C|k|^{-s} \Bigr\} \ .
\end{align}

We will show that if $K$ is chosen large enough, any point starting in $\Omega(K)$ cannot leave
$\Omega(K)$, because the vector field along the boundary $\partial \Omega(K)$ is pointing inward, i.e.,
$\Omega(K)$ is a trapping region. Since the initial data begins in $\Omega(K)$, proving this would prove
the theorem.

We pick a point on $\partial \Omega(K)$ where $\Re(u_i^{(\kb)})$ or $\Im(u_i^{(\kb)})= \pm C|\kb|^{-s}$
for some $\kb \in \mathcal Z$ such that $|\kb|>K$ and some $i \in \{1,2,3\}$. (For definiteness, we shall
assume that $\Re(u_i^{(\kb)})= C|\kb|^{-s}$, but the same line of argument which follows also applies to
the other possibilities.) Then the following inequalities hold when $K$ is chosen large enough:

\begin{align}
\label{eq:ineq} \Bigl|\sum_{\substack{q+r=\kb\\
        q, r \in \mathcal{Z}}}\sum_{\substack{j=1,2,3\\l=1,2,3}}2 \pi
 (\delta_{il}-\frac{\kb_i \kb_l}{|\kb|^2})
   \kb_j \Im(u_l^{(q)}u_j^{(r)})\Bigr|\leq
\sum_{\substack{q+r=\kb\\q, r \in \mathcal{Z}}} \sum_{\substack{j=1,2,3\\l=1,2,3}} 4\pi |\kb_j|
|u_l^{(q)}||u_j^{(r)}| \leq \notag \\
 \sum_{\substack{j=1,2,3\\l=1,2,3}} 4\pi |\kb_j|\Bigl(\sum_{\substack{q \in \mathcal{Z}}}|u_l^{(q)}|^2
 \Bigr)^{1/2}\Bigl(\sum_{\substack{r \in \mathcal{Z}}} |u_j^{(r)}|^2\Bigr)^{1/2} \leq
 \sum_{\substack{j=1,2,3\\l=1,2,3}} 4\pi |\kb_j|E < 4\pi^2\nu|\kb|^\alpha
 \frac{C}{|\kb|^s}=4\pi^2\nu|\kb|^\alpha |\Re(u_i^{(\kb)})|.
\end{align}

This establishes that the vector field points inward along the boundary of $\Omega(K)$ for all $t>0$. So
the trajectory never at any time leaves $\Omega(K)$. Then we have the desired estimate that $|u^{(k)}(t)|
\leq C |k|^{-s}$ for all $t>0$.\qed

\bigskip Just as in the 1999 paper by Mattingly and Sinai \cite{b:MatSi99}, an existence and uniqueness theorem
for solutions follows from our theorem by standard considerations (see
\cite{b:CoFo88,b:DoGi95,b:Temam79}). The line of argument is as follows: By the Sobolev embedding
theorem, the Galerkin approximations are trapped in a compact subset of $L^2$ of the 3-torus. This 
guarantees the existence of a limit point which can be shown to satisfy \eqref{eq:final}, 
where $\mathcal{Z}=\mathbb{Z}^3$. Using the
regularity inherited from the Galerkin approximations, one then shows that there exists a unique
solution to the generalized 3D Navier-Stokes equations where $\alpha > 2.5$.

The inequality \eqref{eq:ineq} in the proof of our Theorem is not necessarily true when $\alpha = 2$. 
Because of this, there is nothing preventing the solutions to \eqref{eq:final} from escaping the region 
$\Omega(K)$ when $\alpha =2$. Hence, there is no logical reason why the standard 3D Navier-Stokes 
equations must always have solutions, 
even when the initial velocity vector field is smooth; if they do always have solutions, it is due to 
probability (see \cite{b:MP02}) and not logic, just like the Collatz $3n+1$ Conjecture and the Riemann 
Hypothesis (see \cite{b:Fei05,b:Fei12}). Of course, it is also possible that there is a counterexample
to the famous unresolved conjecture that the Navier-Stokes equations always have solutions when the initial 
velocity vector field is smooth. But as far as 
the author knows, nobody has ever found such a counterexample.

\end{document}